\documentstyle[]{article}
\def\picture#1by#2(#3){
\vbox to #2 {
  \hrule width #1 height 0pt depth 0pt \vfill \special{picture #3}}
}

\def\scaledpicture#1by#2(#3scaled#4){{
\dimen0=#1  \dimen1=#2
\divide\dimen0 by 1000 \multiply\dimen0 by #4
\divide\dimen1 by 1000 \multiply\dimen1 by #4
\picture \dimen0 by \dimen1 (#3 scaled #4)}}
\def\dfigure#1by#2(#3scaled#4offset#5:#6)
  {\medskip
   \vglue 2mm minus 2mm
   $$
     \hbox{
       \hglue#5
       {\scaledpicture #1 by #2 (#3 scaled #4)}
     }
   $$
   \par\goodbreak
   \vglue 2mm minus 2mm
   \medskip}

\def\qmod#1#2{{\hbox{}^{\displaystyle{#1}}}\!\big/\!\hbox{}_{
\displaystyle{#2}}}

\newfam\msbfam

\font\tenmsb=msbm10
\font\sevenmsb=msbm10 at 7pt
\font\fivemsb=msbm10 at 5pt

\textfont\msbfam=\tenmsb
\scriptfont\msbfam=\sevenmsb
\scriptscriptfont\msbfam=\fivemsb

\def\Bbb{\fam\msbfam\tenmsb}

\def\C{{\Bbb C}}

\def\N{{\Bbb N}}

\def\R{{\Bbb R}}

\def\qed {\hfill\vrule height6pt width6pt depth0pt \bigskip}
\def\map{\longrightarrow}
\def\textmap#1{\mathop{\vbox{\ialign{
                                ##\crcr
    ${\scriptstyle\hfil\;\;#1\;\;\hfil}$\crcr
    \noalign{\kern-1pt\nointerlineskip}
    \rightarrowfill\crcr}}\;}}

\def\textlmap#1{\mathop{\vbox{\ialign{
                                ##\crcr
    ${\scriptstyle\hfil\;\;#1\;\;\hfil}$\crcr
    \noalign{\kern-1pt\nointerlineskip}
    \leftarrowfill\crcr}}\;}}

\newfam\meuffam

\font\tenmeuf=eufm10
\font\sevenmeuf=eufm7

\textfont\meuffam=\tenmeuf
\scriptfont\meuffam=\tenmeuf
\scriptscriptfont\meuffam=\sevenmeuf

\def\germ{\fam\meuffam\tenmeuf}
\def\g{{\germ g}}

\def\kg{{\germ k}}
\def\lg{{\germ l}}
\def\mg{{\germ m}}

\def\tg{{\germ t}}
\def\ug{{\germ u}}

\def\zg{{\germ z}}
\begin{document}
\def\Pr{{\rm Pr}}
\def\tr{{\rm Tr}}
\def\End{{\rm End}}
\def\Aut{{\rm Aut}}
\def\Spin{{\rm Spin}}
\def\U{{\rm U}}
\def\SU{{\rm SU}}
\def\SO{{\rm SO}}
\def\PU{{\rm PU}}
\def\GL{{\rm GL}}
\def\spin{{\rm spin}}
\def\u{{\rm u}}
\def\su{{\rm su}}
\def\so{{\rm so}}
\def\ub{\underbar}
\def\pu{{\rm pu}}
\def\Pic{{\rm Pic}}
\def\Iso{{\rm Iso}}
\def\NS{{\rm NS}}
\def\deg{{\rm deg}}
\def\Hom{{\rm Hom}}
\def\Aut{{\rm Aut}}
\def\h{{\germ h}}
\def\Herm{{\rm Herm}}
\def\Vol{{\rm Vol}}
\def\pf{{\bf Proof: }}
\def\id{{\rm id}}
\def\i{{\germ i}}
\def\im{{\rm im}}
\def\rk{{\rm rk}}
\def\ad{{\rm ad}}
\def\h{{\bf H}}
\def\coker{{\rm coker}}
\def\dv{\bar\partial}
\def\Ad{{\rm Ad}}
\def\RSU{\R SU}
\def\ad{{\rm ad}}
\def\dva{\bar\partial_A}
\def\da{\partial_A}
\def\p{\partial\bar\partial}
\def\sp{\Sigma^{+}}
\def\sm{\Sigma^{-}}
\def\spm{\Sigma^{\pm}}
\def\smp{\Sigma^{\mp}}
\def\oo{{\scriptstyle{\cal O}}}
\def\ooo{{\scriptscriptstyle{\cal O}}}
\def\sw{Seiberg-Witten }
\def\pa{\partial_A\bar\partial_A}
\def\Dr{{\raisebox{0.15ex}{$\not$}}{\hskip -1pt {D}}}
\def\gr{{\scriptscriptstyle|}\hskip -4pt{\g}}
\def\subsetint{{\  {\subset}\hskip -2.45mm{\raisebox{.28ex}
{$\scriptscriptstyle\subset$}}\ }}
\def\nr{\parallel}
\def\ra{\rightarrow}
\newtheorem{sz}{Satz}[section]
\newtheorem{thry}[sz]{Theorem}
\newtheorem{pr}[sz]{Proposition}
\newtheorem{re}[sz]{Remark}
\newtheorem{co}[sz]{Corollary}
\newtheorem{dt}[sz]{Definition}
\newtheorem{lm}[sz]{Lemma}
\newtheorem{cl}[sz]{Claim}

\title{Symplectic stability, analytic stability in non-algebraic complex geometry}
\author{Andrei Teleman\thanks{The starting point of this  article was a research project on ``The
universal Kobayashi-Hitchin correspondence and its applications`` which began in Z\"urich in
  collaboration with Martin L\" ubke and Christian Okonek.}
\date{September 24, 2003}}
\maketitle
\setcounter{section}{-1}

\begin{abstract}
  We give a systematic presentation of the stability theory in the
non-algebraic K\"ahlerian geometry. We introduce the concept of "energy complete
Hamiltonian action". To an energy complete Hamiltonian action of a reductive
group G on a complex manifold one can associate a G-equivariant maximal weight
function and prove a Hilbert criterion for semistability. In other words, for
such actions, the symplectic semistability and analytic semistability
conditions are equivalent.
\end{abstract}

\section{Introduction}

 The factorization problem for group actions in both algebraic
geometry and complex geometry is a very interesting and important subject.  It is
well known that one should impose certain restrictions on the action in order to
get a quotient with good properties. First, following the principles of
the classical theory of invariants, as developed by Mumford, we will
only consider actions of complex reductive groups. Second, it is well
known that, in order to obtain a Hausdorff quotient with an induced complex space
structure satisfying the natural universal property, one should not try to
factorize the whole manifold on which the given reductive group $G$ acts, but
only a certain open part of it (the so called   {\it  stable} locus).  In the
algebraic geometric framework, the  stability condition depends on the
choice of a {\it linearization} of the action in an ample line bundle [MFK], which
is a purely algebraic geometric concept.  Therefore, it is not clear at all how to
generalize   this theory to the K\"ahlerian non-algebraic framework.

The first aim of this paper is to give a systematic presentation of the
different notions of stability in complex K\"ahlerian non-algebraic 
geometry, and to  explain the relations between these notions. There are
two important motivations  for writing this article:
\begin{enumerate} 
\item In the mathematical literature one can find two distinct
stability theories for actions of reductive groups $G$ on K\"ahlerian
manifolds: the symplectic (Hamiltonian) stability and the
analytic stability. In the former theory ([HH],  [Ki])  stability is checked using the
position of the $G$-orbit with respect to the vanishing locus of a moment map
with respect to a maximal compact subgroup of $G$. In  the latter (see for instance
[Mu]), one uses a {\it numerical} criterion, which can be regarded as a
K\"ahlerian version of the Hilbert Criterion in GIT. A well-known comparison
result   states that these two conditions are in fact equivalent. However, there
is no analogous comparison result for the corresponding {\it semi}stability
conditions, and this seems to be  a very delicate point.  Moreover, it is not
clear at all whether the analytical semistability condition is invariant under the
$G$-action\footnote{The difficulty comes from the fact that, in general, one
has no $G$-equivariance property for the maximal weight function
$(x,s)\mapsto  \lambda^s(x)$ (see section 2).} or whether this condition is an
open condition. In particular, one cannot state   that the
analytically-semistable locus has a good quotient.

In fact, the two semistability conditions {\it cannot} be equivalent for
general  Hamiltonian actions on (non-compact) K\"ahler manifolds; one
certainly  needs a {\it completeness condition} for the action.

Indeed, suppose that in a compact manifold $X$ endowed with a Hamiltonian
action with moment map $\mu$, the orbit $Gx_{x_{0}}$ is closed,   $x_0\in \overline{Gx}\cap
\mu^{-1}(0)$, but
$x_0\not\in Gx$. Then, in the open manifold $X\setminus Gx_0$ the point
$x$ will be analy\-tically semistable but no longer symplectically semistable.
\item The analytic   semistability condition, as defined in the
literature, does not have a purely complex geometric character; it depends on the
choice of a   maximal compact subgroup $K$ of the given reductive group $G$ and
on a moment map for the induced $K$-action.  It is not clear  at all that
changing $K$ (and modifying the K\"ahler metric and the  moment map
accordingly) will give the same analytic semistability condition. 

On the other hand, the Hilbert criterion in the algebraic geometric  GIT has
obviously a purely complex geometric character; no differential geometric data are
necessary.
\end{enumerate}

Therefore, the main goals of this article are:

\begin{itemize}
\item To show that, for a large class of actions (which includes all
Hamiltonian actions on compact  K\"ahlerian manifolds and all linear
representations), the analytic semistability is $G$-invariant, has a purely
complex geometric character and is an open condition. Therefore, one can state
 that, for this class of actions, the analytic semistability condition is the
natural extension (to the K\"ahlerian framework) of the algebraic geometric
semistability condition provided by the Hilbert criterion in the algebraic
geometric GIT. 
\item To prove comparison results for our class of Hamiltonian actions relating
the Hamiltonian (semi)stability to the analytic (semi)stability, and identifying
the corresponding quotients.  
\end{itemize}

Via this correspondence we will study carefully the "polystable" orbits, i.
e. the complex orbits which intersect the vanishing locus of the moment map,
and we characterize these orbits too with a numerical (analytic) criterion. Note
that the space of polystable orbits can be identified with the underlying
topological space of the Hamiltonian quotient, hence the polystable orbits are in
fact those which effectively "contribute"  to this quotient.

In the first chapter we review briefly the  main  results of
Heinzner-Huckleberry-Loose [H], [HH], [HHL]  concerning the existence of the Hamiltonian
quotient of a K\"ahler manifold endowed with a Hamiltonian action, and we recall
the first numerical criterion for symplectic stability. Many of the results in
this chapter are well known, but we included short proofs for completeness. In the
second chapter we introduce the concept of {\it energy complete Hamiltonian
action}, concept which plays a fundamental role in  our results. Any Hamiltonian
action on a compact complex manifold and any "linear" Hamiltonian action is
energy complete.  We study the analytic stability, semistability and polystability
conditions with respect to such an action, and we prove the fundamental properties
of the (poly, semi)stable points. The third chapter is dedicated to   comparison
results and to explicit (poly, semi)stability criteria.   

We believe that the energy completeness condition gives the natural framework for
the stability theory in non-algebraic complex geometry.  In the joint paper with L.
Bruasse [BT], we showed that the theory of optimal destabilizing one-parameter
subgroups and a very general Harder-Narasimhan type theorem  can be extended from
GIT to this very large class of holomorphic actions on complex manifolds.

The ideas and the methods of this article can be extended in the infinite dimensional gauge
theoretical framework. This  direction will  be developed in a forthcoming article

\section{Symplectic  stability  and Hamiltonian
K\"ahler quotients}

\subsection{Hamiltonian K\"ahler quotients}

The symplectic (semi)stability condition and the theory of symplectic (Hamiltonian)
K\"ahle\-rian quotients have their roots in the Marsden-Weinstein theory of {\it
symplectic quotients.} The great advantage of this approach is the generality: by the
 results of Heinzner-Huckleberry-Loose, the semistable locus of {\it
any} holomorphic Hamiltonian action on {\it any} K\"ahlerian manifold admits a
{\it good quotient} which can be identified as a topological space with the
corresponding (possibly singular) symplectic quotient. {\it There is no
compactness or completeness condition needed.}

Let $\alpha:G\times X\ra X$ be a holomorphic
action of a reductive group
$G$ on a connected  complex manifold $X$. Suppose that there exists a K\"ahlerian
metric $g$ on $X$ and a
  a maximal compact subgroup $K$ of
$G$ acting by isometries on $(X,g)$. Suppose that the restricted
symplectic action
$\alpha|_{K\times X}:{K\times X}\ra X$  on the symplectic manifold
$(X,\omega_g)$ is {\it Hamiltonian}, i. e. it has a moment map $\mu:X\ra
\kg^\vee$. 
We will suppose for simplicity that $\ker(G\ra\Aut(X))$ is discrete (otherwise
one can factorize $G$ by the connected component $\ker(G\ra\Aut(X))_e$, which is
a reductive normal subgroup of $G$). 

For a point $x\in X$ we denote by $G_x$ (respectively $K_x$) its stabilizer
subgroup with respect to the $G$-action ($K$-action) and by
$\g_x$ (respectively $\kg_x$) its Lie algebra, which is
$$\g_x:=\{u\in\g|\ u^\#_x=0\}\ (\kg_x:=\{u\in\kg|\ u^\#_x=0\})\ .
$$

\begin{dt}\label{SympStab}
A point $x\in X$ is called
\vspace{-2mm}
\begin{enumerate}
\item symplectically $\mu$-{\it stable} if $\g_x=0$ and $Gx\cap
\mu^{-1}(0)\ne\emptyset$.
\item symplectically $\mu$-{\it semistable} if $\overline{Gx}\cap
\mu^{-1}(0)\ne\emptyset$.
\item symplectically $\mu$-{\it polystable} if $Gx\cap
\mu^{-1}(0)\ne\emptyset$.
\end{enumerate}
\end{dt}
We denote by $X^{s}_\mu$, $X^{ss}_\mu$, $X^{ps}_\mu$ the loci of
symplectically $\mu$-stable (respectively semistable, polystable) points.
These loci are obviously
$G$-invariant, so one can speak about stable (semistable, polystable)
$G$-orbits.\\

We refer to   [HH], Lemma 2.4.8, p. 325 for the following important result.
\begin{lm} \label{HH} Let   $H$ be a reductive subgroup of $G$, $L$ a compact
maximal subgroup of $H$ and $K$ a maximal compact subgroup of $G$ which contains
$L$.

Let $\mg$ be the orthogonal complement of $\lg$ in $\kg$ with respect to an
$\ad$-invariant inner product on $\kg$.

Then $K\exp(i\mg)$ is a closed submanifold of $G$ which intersects every
right $H$-congruence class  $\chi\in G/H$ along a unique right $L$-congruence
class
$\lambda\subset\chi$.

Therefore every element $g\in G$ can be decomposed as
$g=k\gamma  h$ with  $k\in K$, $\gamma\in\exp(i\mg)$,  and $h\in H$.
\end{lm}

We will also need the following  well known lemma (see [HH], section 3.2, p. 331). We
include a self-contained proof for completeness.
\begin{pr} \label{RedStab} If $\mu(x)=0$, then $G_x=K_x^\C$.  In particular, 
the stabilizer $G_x$ of any  symplectically polystable point $x$ is a reductive
subgroup of $G$.
\end{pr}
\pf   Suppose that $\mu(x)=0$.

We will show first the infinitesimal version of the the claimed formula, which is
\begin{equation}\label{InfStabComplx}
\mu(x)=0\Rightarrow \g_x=\kg_x^\C
\end{equation}

Let
$w=u+iv\in\g_x$ with $u$, $v\in\kg$.  Then
$$0=\langle w^\#_x,w^\#_x\rangle=\langle u^\#_x,u^\#_x\rangle+   \langle
v^\#_x,v^\#_x\rangle+2\omega( u^\#_x,v^\#_x)=$$
$$\| u^\#_x\|^2+\|
v^\#_x\|^2+2(\iota_{u^\#}\omega)(v^\#_x)
=\| u^\#_x\|^2+\|
v^\#_x\|^2+2  v^\#_x  (\mu^u)\ .
$$
But,
$$v^\#_x  (\mu^u)=\frac{d}{dt}|_{t=0}(\mu^u(\exp(tv)x))=\frac{d}{dt}|_{t=0}
(\mu^{\ad_{\exp(tv)}u} (x))=\mu^{[v,u]} (x)=0\ .
$$
Therefore, the above formula gives     $u^\#_x=v^\#_x=0$, hence
$w\in\kg_x^\C$.

Now we come back to the proof of the equality $G_x=K_x^\C$. The inclusion
$G_x\supset K_x^\C$ is obvious, so let
$g\in G_x$. By Lemma  \ref{HH}, we may decompose $g$ as $g=k\gamma h$ with $k\in
K$,
$\gamma\in \exp(i\mg)$  and $h\in K_x^\C$, where $\mg$ is the orthogonal
complement of
$\kg_x$ in $\kg$ with respect to an $\ad$-invariant inner product.
We get
$$0=\mu(x)=\mu(gx)=\ad_k (\mu(\gamma x))\ ,
$$
hence $\mu(\gamma x)=0$. Write $\gamma=\exp(m)$ with $m\in i\mg$ and consider
the real function $t\mapsto\lambda^m_x(t)=\mu^{-im}(\exp(tm) x)$.

This map is monotone increasing, because, by the properties of the moment map,
$$\frac{d}{dt}(\mu^{-im}(\exp(tm) x))=d\mu^{-im}(J(-im)^\#_{(\exp(tm)
x)})=
$$
$$=\omega_g({-im}^\#_{(\exp(tm) x)}, J({-im})^\#_{(\exp(tm) x)})=g(m^\#_{(\exp(tm)
x)},m^\#_{(\exp(tm) x)})\ .
$$
 
Since $\lambda^m_x(0)=\lambda^m_x(1)=0$, we must   have
$$\frac{d}{dt}|_{t=0} (\lambda^m_x(t))=\nr m^\#(x)\nr^2=0\ ,
$$
hence $m\in\g_x$. But we know that,
since
$\mu(x)=0$, one has $\g_x=\kg_x^\C$. 
Therefore $m\in i\mg\cap \kg_x^\C=\{0\}$, so
$g=kh$ with $k\in K$ and $h\in K_x^\C\subset G_x$. This implies $k\in K\cap G_x=K_x$, so indeed
$g\in K_x^\C$.

\qed

The results of Heinzner-Huckleberry-Loose [H], [HHL], [HH] show that
\begin{thry}\label{KQuot}\hfill{\break}
\vspace{-6mm}
\begin{enumerate}
\item The subsets $X^{s}_\mu$, $X^{ss}_\mu$ are open in $X$.
\item The closure in $X^{ss}_\mu$ of every $\mu$-semistable $G$-orbit contains
a unique $\mu$-polystable orbit. 
\item There is a good quotient $q_\mu: X^{ss}_\mu\ra Q$ with the properties
\begin{enumerate}
\item The induced morphism $ {\mu^{-1}(0)}/{K}\ra Q$ is a homeomorphism,
\item Two $\mu$-semistable $G$-orbits have the same image in $Q$ if and only if
their closures in $X^{ss}_\mu$ are not disjoint, and this happens if and only
if the $\mu$-polystable orbits in their closures in $X^{ss}_\mu$ coincide.
\end{enumerate}
\end{enumerate}
\end{thry}

Therefore, the quotient $Q$ can be identified with the space of polystable
orbits, but, in general, the polystable locus is in general neither open nor closed.
Moreover,   in the algebraic geometric   framework the stable and semistable loci are
both Zariski open. The polystability condition does not appear at
all in the classical GIT. The algebraic geometric version of the polystability
condition was first introduced  in [OST], by generalizing in a natural way the well-known
polystability condition for holomorphic vector bundles.

As formulated above, the concepts of Hamiltonian (semi-, poly-) stability
depends on the choice of three differential geometric objects:
\begin{itemize}
\item a maximal compact subgroup $K$ of $G$.
\item a $K$-invariant K\"ahler metric $g$ on the complex manifold $X$.
\item a moment map $\mu$ for the symplectic action of $K$ on $(X,\omega_g)$.
\end{itemize}

It is very useful to notice that only the class of the triple $(K,g,\mu)$ with
respect to the natural action of $G$ on the set of such triples is essential.
This motivates the following
\begin{dt}\label{DefSymp}  Let $\alpha:G\times X\ra X$ a holomorphic action
of the reductive group $G$ on the complex manifold $X$.  A
\ub{symplectization} of
$\alpha$ is an equivalence class of triples $(K,g,\mu)$ consisting of a 
maximal compact subgroup $K$ of $G$, a $K$-invariant  K\"ahler metric
$g$ on $X$,  and a moment map $\mu$ for
the $K$-action on the symplectic manifold $(X,\omega_g)$.

Two such triples are considered equivalent if there exists $\gamma\in G$ such
that 
$$K'=\Ad_\gamma(K)\ ,\
g'=(\gamma^{-1})^* g\ ,\ \mu'= \ad_{\gamma^{-1}}^t\circ \mu\circ
\gamma^{-1}
$$
\end{dt}
{\bf Remark:} The concept  "{\it symplectization of a holomorphic action}"
should be regarded as the complex geometric analogous   of the algebraic
geometric  concept "{\it linearization of a regular action in an ample line
bundle}" in classical GIT.

It is  convenient to fix a symplectization $\sigma$ of our holomorphic
action $\alpha$, rather than a representative  $(K,g,\mu)$,    and to use the
terminology
"{\it symplectically $\sigma$-stable" ("$\sigma$-semistable",
"$\sigma$-polystable")}, instead of
"symplectically $\mu$-stable" (respectively  "$\mu$-semistable",
"$\mu$-polystable")  for a point
$x\in X$ which satisfies the corresponding condition in Definition \ref{SympStab}.  We
will use the notations $X^{s}_\sigma$, $X^{ss}_\sigma$, $X^{ps}_\sigma$ for the
corresponding subsets of $X$.

 The quotient
$Q$ given by Theorem \ref{KQuot} will be called the {\it K\"ahlerian quotient  of}
$X$ with respect to the symplectization $\sigma$ and will be denoted by
$Q_\sigma$.   Fixing a representative $(K,g,\mu)\in\sigma$ gives a K\"ahler
metric on the smooth part of the quotient $Q_\sigma$.\\

\subsection{A numerical criterion for symplectic stability and poly\-stability}

Let $\alpha:G\times X\ra X$ be a holomorphic action on a complex manifold $X$
and let   $(K,g,\mu)$ be a triple consisting of  maximal compact subgroup $K$
of $G$, a $K$-invariant K\"ahlerian metric $g$ on $X$, and a moment map $\mu$
for the $K$-action on the symplectic manifold $(X,\omega_g)$. 

The analytical
stability condition involves a very important numerical invariant for a system
of data as above, called (for   historical reasons) the maximal weight
function.

For pair $(x,s)$ with $x\in X$ and $s\in i\kg$, consider the path
$c_x^s:[0,\infty)\ra X$  and the map
$\lambda^s_x:[0,\infty)\ra\R$ defined by 
$$c_x^s(t):=e^{ts} x\ ,\ \lambda^s_x(t):=\mu^{-is}(e^{ts} x)\ ,$$ 
where, in general, for $\xi\in\kg$, we use the notation $\mu^\xi$ for the map
$X\ra \R$ given by $y\mapsto \langle \mu(y), \xi\rangle$.
One has (see the proof of Proposition \ref{RedStab})
\begin{equation}\label{derivative}
\frac{d}{dt}\lambda^s_x(t)=g(\frac{d}{dt}c^s_x,\frac{d}{dt}c^s_x)
=g(s^\#\circ c^s_x,s^\#\circ c^s_x)\ , 
\end{equation}
hence, for fixed $s\in i\kg$, the map $\R\times X\ni(t,x)\mapsto
\lambda^s_x(t)\in\R$ is increasing with respect to the first argument. We put
$$\lambda^s(x):=\lim_{t\ra\infty}\lambda^s_x(t)\in\R\cup\{\infty\}\ ,
$$
and we call $\lambda^s_x$ the {\it maximal weight} of $x$ in the direction $s$. 
The above formula (\ref{derivative}) shows that
\begin{equation}\label{energy}
\lambda^s(x)=\lambda^s_x(0)+E_g(c^s_x)\ ,  
\end{equation}
where $E_g$ stands for the energy with respect to the metric $g$.

We will need the following important notion: 
\begin{dt}\label{Psi} The {\it integral of the moment
map}  is the map $\Psi:X\times G\ra\R$ satisfying the conditions
\begin{enumerate}
\item $\Psi(x,e)=0$ for all $x\in X$.
\item $\Psi$ is $K$-invariant from the left, i. e. $\Psi(x,kg)=\Psi(x,g)$ for
all $x\in X$, $g\in g$, $k\in K$.
\item $\Psi(x,gh)=\Psi(x,h)+\Psi(hx,g)$ for all $x\in X$, $g,\ h\in G$.
\item $\frac{d}{dt}\Psi(x,e^{ts})=\mu^{-is}(e^{ts}x)$ ($=\lambda^s_x(t))$. 
\end{enumerate}
\end{dt}
Using the fact that the map $K\times i\kg\ra G$ given by $(k,s)\mapsto k e^s$
is a diffeomorphism (see  for instance [HH]), it is not difficult to see that these conditions
determine a unique function $\Psi:X\times G\ra\R$.
By (\ref{derivative})  and   the  
properties listed in Definition \ref{Psi}, one notes that
\begin{re}{\ }\vspace{-2mm} 
\label{PsiProp} 
\begin{enumerate}
\item $\Psi|_{X\times K}\equiv 0$.
\item  For any fixed 
$x\in X$, $s\in i\kg$, the real function $t\mapsto \Psi(x,e^{st})$ is convex.
\item For any $x\in X$ the restriction $\Psi|_{G_x}$ is an $\R$-valued group
morphism.
\end{enumerate} 
\end{re}
Using the  properties 3, 4 in Definition \ref{Psi}, one gets immediately the
following simple but important
\begin{re}\label{Crit} Let $x\in X$. The
following conditions are equivalent:
\begin{enumerate}
\item $g_0$ is a critical point of the map $\Psi(x,\cdot)$.
\item $\mu(g_0x)=0$
\end{enumerate}
\end{re}

We will need the following well known  lemma (see [Mu]). 
We include a proof for completeness.

\begin{lm}\label{LinProp} Fix   an $\ad_K$-invariant metric on $\kg$ and a  
subspace
$V\subset i\kg$. For a point $x\in
X$ the following conditions are equivalent
\begin{enumerate} 

\item The map $\Psi(x,\exp(\cdot))$ is linearly proper on 
 $V$, i. e.
there exist positive constants $c_1$, $c_2$ such that
$$\nr s \nr\leq c_1\Psi(x,e^s)+c_2\ ,\ \forall s\in V\ .$$

\item $\lambda^s(x)>0$ for all $s\in V\setminus\{0\}$.
\end{enumerate}
\end{lm}
\pf 1. $\Rightarrow$ 2.: The inequality in 1. gives for any  $s\in
V$, $t\in\R$
$$t\nr s \nr\leq c_1\Psi(x,e^{st})+c_2\ .
$$
This shows that
$$\frac{d}{dt}\Psi(x,e^{st})=\lambda^s_x(t)>0
$$
for $s\ne 0$ and sufficiently large $t\in\R$, hence $\lambda^s(x)>0$.\\
2.$\Rightarrow 1$. Suppose that there didn't exist any positive constants
$(c_1,c_2)$ with the required property. It would follow that there exist a
sequence of $(s_n)_n$ in
$V$
$$\nr s_n \nr> n\Psi(x,e^{s_n})+n^2\ .$$
Note that
$$\lim_{n\ra\infty} \nr s_n\nr=\infty\ ,$$
because, if not,  $(s_n)_n$ would have a bounded subsequence
$(s_{n_m})_{m}$. But then
$(\Psi(x,e^{s_{n_m}}))_m$ would be also bounded, and this obviously contradicts
the above inequality.

We get
$$\frac{\Psi(x,e^{s_n})}{\|s_n\|}<\frac{1}{n}
$$
Put $l_n:=\nr s_n\nr$, $u_n:=\frac{s_n}{\nr s_n\nr}$, and choose
$t_0\in\R$. The
convexity property of the function $\Psi$ (Remark \ref{PsiProp}) gives
$$\Psi(x,e^{l u})\geq \Psi(x,e^{t_0 u})+(l-t_0)\lambda_{x}^u(t_0)
$$
for every $u\in i\kg$,  $l\geq t_0$.  We obtain
$$\frac{\Psi(x,e^{t_0 u_n})+(l_n-t_0)
\lambda_{x}^{u_n}(t_0)}{l_n}\leq\frac{\Psi(x,e^{s_n})}{\|s_n\|}<\frac{1}{n}\ .
$$
The sequence $(u_n)_n$ has a subsequence which converges to, say
$u_0\in V$ which
must have $\nr u_0\nr=1$. Taking the limit of the right hand term, we get
$\lambda^{u_0}_{x}(t_0)\leq 0$. But this implies
$\lambda^{u_0}(x)\leq 0$, which
contradicts the hypothesis.
\qed

The following simple lemma will play a crucial role in the next chapter.

\begin{lm}\label{Dec}  Let $K$ be a maximal compact subgroup of $G$ and let $g\in
G$,
$s\in\kg$ such that $\ad_g(s)\in\kg$.  Decompose $g$  as $g=kh$, where
$k\in K$ and $h\in \exp(i\kg)$. Then $\ad_h(s)=s$.
\end{lm}
\pf  Since  $\ad_g(s)\in\kg$, one
has $\sigma:=\ad_h(s)=\ad_{k^{-1}}(\ad_g(s))\in\kg$.
Choose an embedding $G\hookrightarrow GL(r,\C)$ mapping $K$ to $U(r)$.
Then the image of $h$ is Hermitian with positive eigenvalues, whereas the
images of
$s$ and
$\sigma$ are anti-Hermitian.  We get
$$-\ad_h(s)=-\sigma=\sigma^*=\ad_h(s)^*=-\ad_{h^{-1}}(s)\ ,
$$
hence $\ad_{h^2}(s)=s$. Therefore the eigenspaces of
$h^2$ (which are the eigenspaces of $h$)  are invariant under $s$, so that
one also has $\ad_h(s)=s$.
\qed

If   $\lg$ is a subset of the Lie algebra $\kg$ of a 
Lie  group $K$, we will denote by $Z_K(\lg)$ 
(respectively $z_\kg(\lg)$) the centralizer of   $\lg$  in $K$
(respectively $\kg$). The Lie algebra of $Z_K(\lg)$ is  $z_\kg(\lg)$.
\begin{re}\label{centralizers}
If $K$ is a maximal compact subgroup of a complex reductive group $G$ and $\lg\subset\kg$, then
$Z_K(\lg)$ is a maximal compact subgroup of the reductive group $Z_G(\lg)$.
\end{re}
\pf One obviously has $z_\g(\lg)=z_\kg(\lg)\otimes\C$, so it suffices to prove
that
$Z_G(\lg)=Z_K(\lg)\exp(i z_\kg(\lg))$.  Let
$g\in Z_G(\lg)$, and decompose $g$ as $g= kh$ with $k\in K$ and $h\in
\exp(i\kg)$. By Lemma \ref{Dec}, it follows that $h\in
Z_G(\lg)$, hence   $k\in Z_K(\lg)$. But, using an embedding $G\ra
GL(r,\C)$ mapping $K$ to $U(r)$, one gets easily that
$Z_G(\lg)\cap\exp(i\kg)=\exp(iz_\kg(\lg))$.
\qed

\begin{re}\label{MomCent}
Let $\mu$ be a moment map for an action of a compact group $K$ on a
symplectic manifold $(M,\omega)$, and let $x\in M$. Then, via an
identification $\kg^\vee=\kg$ given by an $\ad$-invariant inner product on
$\kg$, one has  $\mu(x)\in z_{\kg}(\kg_x)$.
\end{re}
\pf The equivariance property of the moment map gives 
$$0=\frac{d}{dt}|_{t=0}\ \mu(\exp(tu)x)=\frac{d}{dt}|_{t=0}\
\ad_{\exp(tu)}(\mu(x))=[u,\mu(x)]\ .
$$
\qed

Now we can prove the following important

\begin{lm}\label{lemma}
Let $\alpha:G\times X\ra X$ an holomorphic action of a complex reductive group
on a complex manifold, $K$ a maximal compact subgroup of $G$, $g$ a
$K$-invariant K\"ahler metric on $X$ and $\mu$ a moment map for the $K$ action.
Let $x\in X$ such that $\lambda^s(x)\geq 0$ for all $s\in iz_{\kg}(\kg_x)$
and $\lambda^s(x)>0$ for all 
$s\in iz_{\kg}(\kg_x) \setminus iz(\kg_x)$. Then 
\begin{enumerate}
\item There exists $s_0$ in the
orthogonal complement $i[z(\kg_x)^{\bot_{z_{\kg}(\kg_x)}}]$  of $iz(\kg_x)$
in $iz_{\kg}(\kg_x)$ with respect to an $\ad$-invariant inner product on
$\kg$, such that
$$\mu(\exp(s_0) x)=0\ .$$
\item  $x$ is symplectically polystable.

\end{enumerate}
\end{lm}
\pf

1. Consider the restricted action $\alpha':=\alpha|_{Z_K(\kg_x)\times X}$ of
the centralizer $Z_K(\kg_x)$ and the induced moment map $\mu'$ for this
action. Note that the   function
$\Psi'$ associated with the  triple $(Z_K(\kg_x),g,\mu')$ is just the
restriction to $X\times Z_G(\kg_x)$  of the function $\Psi$ corresponding to
$(K,g,\mu)$ (see Remark \ref{centralizers}).

 Apply Lemma \ref{LinProp} to the triple $(Z_K(\kg_x),g,\mu')$ taking
$V:=i[z(\kg_x)^{\bot_{z_{\kg}(\kg_x)}}]$.   It follows that
$\Psi'(x,\exp(\cdot))$ is linearly proper on this space, hence there exist 
positive constants
$c_1$, $c_2$ such that

\begin{equation}\label{ineq}
\nr s\nr\leq c_1\Psi'(x,e^s)+c_2\ ,
\end{equation}
for all $s\in V$.  This inequality implies that
$\Psi'(x,\exp(\cdot))$ is bounded  from below on $V$.  Put
$$m:=\inf_{s\in V}\Psi'(x,\exp(\cdot)) \ .$$

Let $(s_n)_n$ be a sequence in $V$ such that
$\Psi'(x,e^{s_n})\ra m$.  By (\ref{ineq}) it follows that $(s_n)_n$ is
bounded, so it has a
   subsequence which converges to, say, $s_0\in V$. One
gets
\begin{equation}\label{s_0}
\Psi'(x,e^{s_0})= m\ .
\end{equation}

We claim that in fact
\begin{equation} \label{m}
m=\inf_{g\in Z_G(\kg_x)}\Psi'(x,g)
\end{equation}

Indeed, by  Lemma \ref{HH} applied to  the reductive subgroup
$H:=[\exp(z(\kg_x)]^\C$ (the complexification of the connected component of
$e$ in $Z(K)$) of $Z_G(\kg_x)$, it follows that  any element  
$g\in Z_G(\kg_x)$ can be written as
$g=k\gamma   h $ with  $k\in Z_K(\kg_x)$,  $\gamma\in\exp(V)$
and $h\in H$. We have
$$\Psi'(x,g)=\Psi'(x,k\gamma h)=\Psi'(x,\gamma
h)=\Psi'(x,h)+\Psi'(hx,\gamma)=
\Psi'(x,\gamma)\ ,
$$
because $\Psi'(x,h)=0$. To see this recall that, by Remark \ref{PsiProp}   the restriction
$\Psi'(x,\cdot)|_{H}$ is an $\R$-valued group morphism which vanishes on
the maximal compact subgroup $\exp(z(\kg_x))=Z(K_x)^e$.  But the derivative
of this morphism in the
$s\in iz(\kg_x)$ direction is
$\mu^s(x)=\lambda^s(x)\geq 0$, by assumption. Therefore
$d_e\Psi'(x,\cdot)|_{iz(\kg_x)}$ is an
$\R_{\geq 0}$-valued real linear form, hence  it vanishes.

 This proves the claimed formula (\ref{m}). From (\ref{s_0}) and (\ref{m}) we get that
$e^{s_0}$ is a critical point of the map
$\Psi'(x,\cdot)$, hence $\mu'(\exp(s_0)x)=0$, by Remark \ref{Crit}.

The point is that, in our situation,   $\mu'(\exp(s_0)x)=0$ implies
the stronger relation $\mu(\exp(s_0)x)=0$. Indeed,  since 
$$G_{\exp(s_0)x}=\exp(s_0) G_x \exp(-s_0)$$
and $\exp(s_0)\in Z_G(\kg_x)$ one gets $ \g_{\exp(s_0)x}\supset \kg_x$, so  
$$\kg_x=\kg_x\cap\kg  \subset \g_{\exp(s_0)x}\cap \kg=\kg_{\exp(s_0)x}\ .$$

Therefore, by Remark \ref{MomCent}, $\mu(\exp(s_0)x)\in z_{\kg}(\kg_x)$.
But, via our identification $ \kg^\vee=\kg$, $\mu'(\exp(s_0)x)$ is just the
 orthogonal projection of $\mu(\exp(s_0)x)$ on $z_{\kg}(\kg_x)$, hence
$\mu (\exp(s_0)x)=\mu'(\exp(s_0)x)=0$.
\\ \\
2. This follows immediately from 1.\\
\qed
\\ 

The following   numerical criterion is well-known  in the stable case (see
[Mu]).  
\begin{pr} \label{NumCrit} Let $x\in X$.  
\begin{enumerate}
\item The following conditions are equivalent:
\begin{enumerate}
\item $x$ is symplectically $\mu$-stable.
\item $\lambda^s_x\geq 0$ for all  $s\in i\kg$ with strict inequality for 
$s\in i\kg\setminus\{0\}$.    
\end{enumerate}
\item  The following conditions are equivalent:
\begin{enumerate}
\item $x$ is symplectically $\mu$-polystable.
\item   There exists $g\in G$ such that $\lambda^s_{gx}\geq 0$ for all 
$s\in i\kg$ with strict inequality for 
$s\in i\kg\setminus i\kg_{gx}$.
\end{enumerate}   
\end{enumerate}
\end{pr}
\pf

$1a)  \Rightarrow 1b)$:   Let $x$ be a symplectically stable point and choose
$g_0\in G$ such that $x_0:=g_0x\in  \mu^{-1}(0)$. We   prove first that $x_0$
has the claimed property, namely that
$\lambda^s(x_0)\geq 0$ for all $s\in i\kg$ with strict inequality for $s\ne 0$.
But
$\lambda^s(x_0)\geq
\lambda^s_0(x_0)=0$, because
$\mu(x_0)=0$.

For
$s\ne 0$ we have
$s\not\in
\g_{x_0}$, so
$s^\#(x_0)\ne 0$, which
implies
$$\lambda^s(x_0)>\lambda^s_{x_0}(0)=0\ .$$
This proves the claimed property for $x_0$. Unfortunately these properties do
not appear to be  $G$-invariant, so one cannot deduce directly that the same is
true for $x$. By Lemma
\ref{LinProp}, the map $\Psi(x_0,\exp(\cdot))$ is linearly proper on   $i\kg$. 
One can write
$$\Psi(x,\exp(s))=\Psi(g_0^{-1}x_0,\exp(s))=
\Psi(x_0,\exp(s)g_0^{-1})-\Psi(x_0,g_0^{-1})\ .
$$
Write $\exp(s)g_0^{-1}=k(s)\exp(v(s))$ with $k(s)\in K$ and $v(s)\in i\kg$.
With these notations one has $\Psi(x_0,\exp(s)g_0^{-1})=\Psi(x_0,\exp(v(s)))$. It
is easy to prove  (see [Mu]) an estimate of the form
$$\| s\| \leq a\| v(s)\| +b\ ,\  \forall s\in i\kg\ .
$$
Therefore the linearly properness of $\Psi(x,\exp(\cdot))$ on
$i\kg$ follows from the linearly properness of $\Psi(x_0,\exp(\cdot))$ on
$i\kg$. Applying again Lemma \ref{LinProp}, we get the desired property for $x$.

$1b)  \Rightarrow 1a)$:  Suppose that $1b)$ holds for $x$. First of all notice
that, by Lemma
\ref{lemma}, there exists $s_0\in i\kg$ such that $\mu(\exp(s_0)x)=0$.  It
remains to show that $\g_x=\{0\}$. Put
$x_0:=\exp(s_0)x$. Arguing as above we see that $\Psi(x,\exp(\cdot))$ and
$\Psi(x_0,\exp(\cdot))$  are both linearly proper on $i\kg$, hence 1b) also
holds for $x_0$. We know by Proposition \ref{RedStab} that $G_{x_0}=K_{x_0}^\C$,
hence it suffices to show that $i\kg_{x_0}=0$.  If $s\in i\kg\setminus\{0\}$, one
has $\lambda^s_{x_0}(0)=\mu^{-is}(x_0)=0$ and $\lambda^s(x_0)>0$, hence the path
$c_{x_0}^s$ cannot be constant, so $s\not\in i\kg_{x_0}$.

$2a)  \Rightarrow 2b)$ 

Let $x_0\in Gx$ such that  $\mu(x_0)=0$. 
The same method as in the case of stability, gives
$$\lambda^s(x_0)=0\ {\rm for}\ s\in i\kg_{x_0}\ ,\ \lambda^s(x_0)>0\ {\rm
for}\ s\in i\kg\setminus\i\kg_{x_0}\ .
$$
$2b)  \Rightarrow 2a)$ 

Put $x':= gx$, where $g\in G$ has the   property   in
$2b)$.

By Lemma \ref{lemma},  there exists $s_0\in i\kg$ such that
$\mu(\exp(s_0)x')=0$. Thus $x$ is symplectically polystable.
\qed 
\vspace{3mm}\\
{\bf Remark:} The numerical criterion   provided by Proposition
\ref{NumCrit} is not satisfactory for the following important  reasons:
\begin{itemize}
\item It depends essentially on the choice of a particular triple $(K,g,\mu)$,
not only on its equivalence class (the symplectization defined by this triple).
Therefore, it does not have a purely complex geometric character.
\item In general, for a polystable point $x$, one might have $\lambda^s(x)=0$
even for vectors $s\not\in i\kg_x$. Therefore, in order to test whether a
 point $x$  is polystable or not, one has to control all      the "maximal
weights"
$\lambda^s_{gx}$ as $g$ vary in $G$, so this is not an intrinsic criterion in
terms of the given point $x$.
\item This criterion does not  provide any numerical characterization  
of symplectic
\ub{semi}stability.
\end{itemize}

One of our main goals is to address all these  issues, and to give   intrinsic,
purely  complex geometric numerical criteria for stability, polystability and
semistability (see Theorems \ref{FirstComp}, \ref{semistability}). In order to
get stronger comparison results,  we will have to assume that the triple
$(K,g,\mu)$ satisfies a certain completeness condition, which we will call energy
completeness. Note that for general Hamiltonian actions, analytic semistability
does not imply symplectic semistability and there is no  way to construct a good
quotient of the analytically semistable locus.

\section{Analytic stability, semistability and polystability}

Analytic stability is a purely numerical condition, so it is very useful for
practical reasons.  The analytic stability condition  is the complex geometric
analogue of the numerical condition   in the Hilbert criterion in
classical GIT.

\subsection{ The function $\lambda$ associated
with an energy complete symplectization}

\paragraph{The cone of Hermitian type vectors.} The set $\Hom(\C^*,G)$ of one
parameter subgroups of $G$ can be identified with a subset $\Lambda(G)$ of $\g$
via the map 
$$\lambda\mapsto  d_1(\lambda)(1)=\frac{d}{dt}\big|_{t=0}(\lambda(e^t))\ .$$

In order to formulate a numerical stability condition in  complex non-algebraic
geometry one needs a larger subset of $\g$ whose elements can be interpreted as
"non-algebraic" one parameter subgroups of $G$.

\begin{dt}\label{HermType} Let $G$ be   complex reductive group. An element $s\in
\g$ will be called of Hermitian type if it satisfies one of the following
equivalent properties:
\begin{enumerate}
\item There exists a complex torus $C\subset G$ such that $s\in i\tg$,
where $\tg$ is the Lie algebra of the (unique) maximal compact subgroup $T$
of $C$.
\item The closure of the real one parameter subgroup of $G$ determined by
$is\in\g$ is compact.

\item There exists a compact subgroup $K\subset G$ such that $s\in i\kg$.
\item There exists $r\in\N$ and an embedding $\rho:G\hookrightarrow GL(r,\C)$
such that
$\rho_*(s)$ is Hermitian.
\item There exists   $r\in\N$ and an embedding $\rho:G\hookrightarrow GL(r,\C)$
such that the matrix $\rho_*(s)$ is diagonalizable and has real eigenvalues.
\item For every embedding $\rho:G\hookrightarrow GL(r,\C)$
the matrix $\rho_*(s)$ is diagonalizable and has real eigenvalues.
\end{enumerate}
\end{dt}

\begin{re}\label{inv} \hfill{\break}
1. The set  $H(G)\subset\g$ of elements of Hermitian type
is obviously  invariant under the adjoint action of $G$ on $\g$.\\
2. If
$s\in H(G)$ then the associated endomorphism
$[s,\cdot]\in\End(\g)$ is diagonalizable and has only real eigenvalues.
\end{re}

  One can associate to every $s\in
H(G)$ a parabolic subgroup
$G(s)\subset G$  in the following way:

We put
$$G(s):=\{g\in G|\ \lim\limits_{t\ra\infty}e^{st}ge^{-st}\ {\rm exists\
in}\ G\}\ .
$$
The group $G(s)$ fits in the exact sequence 
$$1\map U(s)\map G(s)\map Z(s)\map 1\ ,
$$
 where $Z(s)$ is the centralizer of $s$ in $G$ and $U(s)$ is the
unipotent subgroup
$$U(s):= \{g\in G|\ \lim\limits_{t\ra\infty}e^{st}ge^{-st}=e\}\ .
$$
Moreover, it  decomposes as a semi-direct product
\begin{equation}\label{SemProd}
G(s)=Z(s)\cdot U(s)\ .
\end{equation} 
The groups $G(s)$, $U(s)$ depend only on the semisimple part $s_0$ of $s$.  The
Lie algebras of
$G(s)$,
$Z(s)$ and $U(s)$ are

$$\g(s):=\bigoplus_{\beta\leq 0} {\rm Eig}([s,\cdot],\beta)\ ,\
\zg(s):=  \ker ([s,\cdot])\ ,\ \ug(s):=\bigoplus_{\beta< 0} {\rm
Eig}([s,\cdot],\beta)\ .
$$
\begin{pr}\label{FreeTran} The adjoint representation defines an affine
representation of
$U(s)$ on the affine subspace $p_{\zg(s)}^{-1}(s)\subset\g(s)$. This
action is free and transitive.
\end{pr}
\pf It  is easy to see that the differential in $e$ of the map
$\varphi_\zeta:U(s)\ra p_{\zg(s)}^{-1}(s)$ given by $u\mapsto
\ad_u(\zeta)$ is invertible, for any $\zeta\in p_{\zg(s)}^{-1}(s)$.  It
follows that $\varphi_\zeta$ is
\'etale. The stabilizer of any point $\zeta\in p_{\zg(s)}^{-1}(s)$ is
trivial, because it is a Zariski closed 0-dimensional subgroup of a
unipotent group. Therefore $\varphi_\zeta$ is an algebraic isomorphism from
$U(s)$ (which, as an algebraic variety,  is isomorphic to an affine
space) onto a Zariski open subset of the affine space
$p_{\zg(s)}^{-1}(s)$. It follows that $\varphi_\zeta$ must be
surjective.

\begin{co}\label{EqRel} Let $s,\ \sigma\in H(G)$. \\
1. The following three conditions are
equivalent:\\
i) $\sigma\in \g(s)$ and $p_{\zg(s)}(\sigma)=s$.\\
ii) $s$ and $\sigma$ are conjugate under the adjoint action of
$U(s)$.\\ iii) $s$ and $\sigma$ are conjugate under the adjoint action
of
$G(s)$.\\
2. If one of these conditions is satisfied then $G(s)=G(\sigma)$. \\
3. The equivalent conditions $i)$ -- $iii)$ define an
equivalence relation $\sim$ on $H(G)$. The equivalence class of $s\in
H(G)$ is the affine subspace $p_{\zg(s)}^{-1}(s)$ of the Lie
algebra $\g(s)$.
\end{co}

\begin{co}\label{Ret} Let $K$ be a maximal compact subgroup of
$G$. Then $i\kg\subset H(G)$ is a complete system of representatives for
$\sim$. Mapping any $s\in H(G)$ to  the representative in $i\kg$ of its
equivalence class, gives a   retraction
$\sigma_K:H(G)\ra\i\kg$, which induces a homeomorphism
$\left[\qmod{H(G)}{\sim}\right]
\rightarrow i\kg$.
\end{co}
\pf Let $K_0$ be a maximal compact subgroup of $G$ such that
$s\in i\kg_0$, and let $g\in G$ such that $\Ad_g(K_0)=K$. Therefore
$\ad_{g}(s)\in i\kg$. Decompose
$g$ as
$g=kb$, where $k\in K$ and $b$ belongs to a Borel subgroup of $G$ contained
in  $G(s)$.

Then
$s\sim\ad_b(s)=\ad_{k^{-1}}\ad_{k}\ad_{b}(s)=\ad_{k^{-1}}(\ad_{g}(s))\in
i\kg$.

Now suppose that $s$, $s'\in i\kg$ and that $s\sim s'$. It follows that
$s'-s\in \ug(s)\cap i\kg$. But $\ug(s)$ is a complex Lie algebra,
hence $i(s'-s)\in \ug(s)\cap  \kg$. This would imply that the closure of
the real 1-parameter subgroup generated by $i(s'-s)$ is compact and
contained in $U(s)$. But an unipotent group contains no compact
subgroups, which completes the proof.
\qed
\begin{co}\label{embedding}
An  embedding
$G\hookrightarrow G'$ of reductive Lie groups  induces   injections
$$H(G)\hookrightarrow H(G')\ ,\ \qmod{H(G)}{\sim_G}\hookrightarrow
\qmod{H(G')}{\sim_{G'}}\ .
$$
\end{co}
\pf  It is easy to see that, when $s\in H(G)$, then $\ug_{G'}(s)\cap\g=\ug_G(s)$,
  where $\ug_G(s)$ stands for the unipotent algebra associated to $s$,
regarded as an element of $H(G)$.  Therefore an element $\sigma\in H(G)$
is equivalent to
$s$  as elements in   $H(G)$ if and only they are equivalent  as
elements in   
$H(G')$.
\qed
\\
{\bf Example:} Consider the case $G=GL(r,\C)$. The data of an equivalence class
 in  $H(G)$ is equivalent to the data of a pair $(\Phi,\eta)$,
where $\Phi$ is a filtration of the form
$\{0\}\subset F_1\subset\dots\subset F_k=\C^r$ and $\eta$ is an
increasing sequence
$\eta_1<,\dots,<\eta_k$ of real numbers. An element $s\in\g l(r,\C)$
belongs to the equivalence class defined by $(\Phi,\eta)$ iff it
is diagonalizable, its spectrum is
$\eta$, and the associated eigenspace filtration $(\oplus_{j=1}^i({\rm
Eig}(s,\eta_j))_i$ of $\C^r$ is
$\Phi$.
\begin{pr}\label{CompactOrbit} Let $G$ be reductive complex group and let $s\in
H(G)$ and let
$K$ be a maximal compact subgroup of $G$ such that $s\in i\kg$. Then the image of
the orbit $\ad_G(s):=\{\ad_g(s)|\ g\in G\}$ of $s$ via the map
$$p_K:H(G)\map \qmod{H(G)}{\sim}\ \stackrel{\simeq}{\map} i\kg
$$
is the compact orbit $\ad_K(s):=\{\ad_k(s)|\ k\in K\}$ of $s$ in $i\kg$.
\end{pr}
\pf  The inclusion of $\ad_K(s)$ in the projection of $\ad_G(s)$ is clear. 
Conversely, let $s'\in p_K(\ad_G(s))$. This means that $s'\in i\kg$ and there
exists $g\in G$ such that $s'\sim \ad_g(s)$. Therefore there exists
$\gamma\in G(\ad_g(s))$ such that $s'=\ad_{\gamma g}(s)$.  Decompose $g':=\gamma
g$ as $g'= k' h'$ with $k'\in K$ and $h'\in\exp(i\kg)$.  By Lemma \ref{Dec}, it
follows that $\ad_{h'}(s)=s$, hence $s'=\ad_{k'}(s)$.
\qed
\\ \\

\paragraph{Energy complete symplectizations.} Let now   $\alpha:G\times X\ra X$
be an action of a complex reductive    group $G$
on a complex manifold
$X$.

To every pair $(s,x)\in  H(G)\times X$ we associate a curve
$c^s_x:[0,\infty)\ra F$ given by $c^s_x(t):= e^{ts} x$.

\begin{dt} \label{EnComp} A symplectization $\sigma$ of the action $\alpha$ will
be  called {\it energy-complete} if, choosing a  representative
$(K,g,\mu)\in\sigma$, the following holds:\vspace{3mm}\\
$$\forall x\in X\ \forall s\in i\kg \ \big( E_g(c^s_x)<\infty\Rightarrow
\hbox{the curve}\  c_x^s \  \hbox{has a limit as}\ t\rightarrow\infty\big)\ .
\eqno{(C)}$$
\end{dt}

Here we denoted by $E_g$ the energy with respect to the metric $g$.
Using the obvious identities
$$c^{\ad_\gamma(s)}_{\gamma x}=\gamma c_x^s\ ,\
E_{(\gamma^{-1})^*g}(c^{\ad_\gamma(s)}_{\gamma x})=E_g(c^s_x)\ ,
$$
one checks that the condition $(C)$ does not depend on the choice of the
representative $(K,g,\mu)$ of $\sigma$.

\begin{pr}\label{EnCompl} \hfill{\break}
1. A symplectization $\sigma$ is energy complete if and only if, choosing  any
representative $(K,g,\mu)\in\sigma$, the following implication holds:
\begin{equation}\label{Cweak}\forall x\in X\ \forall s\in i\kg\ \left(
E_g(c^s_x)<\infty\Rightarrow
\exists (t_n)_{n }\ \ra\infty  \hbox{ s. t. }
 c^s_x(t_n)_{n } \hbox{ converges}\right)  .
\end{equation}
2. Any symplectization of an action
$\alpha:G\times X\ra X$ of a complex reductive group on a compact complex
manifold is energy complete.\\ 3. Let $\alpha:G\times V\ra V$ be a linear action
of a complex reductive group on a finite dimensional Hermitian   vector
space $(V,h)$ and let $K$ be  a maximal compact subgroup of $G$ which leaves $h$
invariant. Then any symplectization   $[(K,h,\mu)]$ of $\alpha$ is
energy complete. 
\end{pr}
\pf

1. It's clear that the condition $(C)$ in the definition implies (\ref{Cweak}).
Conversely, suppose that $(C)$ holds,   let $(x,s)\in X\times i\kg$ and let
$(t_n)_n$ be a sequence  in
$\R_{\geq 0}$ with $t_n\ra\infty$ such that $c_x^s(t_n)$ converges to a point,
say $x_0$, in $ X$. We will show that
\begin{equation}\label{limit}
\lim\limits_{t\ra\infty} c_x^s(t)=x_0\ .
\end{equation}

Consider the compact torus $T=\overline{\{\exp(its)|\ t\in\R\}}\subset K$, and
let $T^\C\subset G$ be its complexification.

Let $\iota_\tg:\tg\hookrightarrow \kg$ be the Lie algebra monomorphism induced
by the inclusion $T\hookrightarrow K$. The maps
$\mu_T=\iota_\tg^*\circ \mu:X\ra\tg^*$, $\mu'_T:=\mu_T-\mu_T(x_0)$ are   moment
maps for the  induced $T$ action on $X$, and obviously $\mu_T'(x_0)=0$.
Using the results in [HH],   Theorem 3.3.14 p. 343 and section 4.1 p. 345 -
346, one can find a
$T^\C$-invariant Stein neighborhood $U$ of $x_0$, a linear representation
$\rho:T^\C\ra GL(V)$ on a finite dimensional vector space $V$, and a
$T^\C$-equivariant isomorphism
$\theta:U\ra W$ on a closed $T^\C$-invariant complex subspace $W$ of $V$.  Let
$R\subset\Hom(\tg^\C,\C)$ be the root set of $\rho$ and
$V=\oplus_{\chi\in R} V\chi$ be the eigenspace decomposition of $V$ with respect
to $\rho$. Since
$U$ is   open, one has $c^s_x(t_n)=\exp(t_n s)x\in U$ for sufficiently large $n$,
hence, since $U$ is also $T^\C$-invariant, one gets $x\in U$. Putting
$v:=\theta(x)$, one can write
\begin{equation}\label{EigenDec}
\theta(c^s_x(t))=\theta(\exp(ts)x)=\rho(exp(ts))\theta(x)=\sum_{\chi\in
R} e^{t\chi(s)}v_\chi\ .
\end{equation}
Since $c^s_x(t_n)\ra x_0$, it follows that $\theta(c^s_x(t_n))\ra
v_0:=\theta(x_0)$  so that, by (9), one must have
$$\chi(s)>0\Rightarrow v_\chi=0\ .
$$
But this implies $\lim_{t\ra\infty}\theta(c^s_x(t))=v_0$,
hence $\lim_{t\ra\infty} c^s_x(t) =x_0$.
\\
\\
2. Follows easily from 1.
\\
\\
3.    The standard moment
map for the
$K$-action on $V$  is 
$$\mu_0(v)=\varsigma^*_\alpha(-\frac{i}{2}v\otimes v^*)\ ,
$$
where $\varsigma_\alpha:\kg\ra u(V)$ is the morphism induced by the
representation $K\ra U(V)$ induced by $\alpha$. Any other moment map has the
form
$$\mu_\tau=\mu_0-i\tau
$$
where $\tau$ is a central element in $i\kg^\vee$.

Let $s\in i\kg$ and $V=\oplus_{j=1}^k V_j$ the decomposition of $V$
as a the direct sum of $\varsigma_\alpha(s)$. Therefore $\varsigma_\alpha(s)
|_{V_j}=s_j\id_{V_j}$, where
$s_j$ are the   eigenvalues of  $\varsigma_\alpha(s)$.

Let $v\in V$. Decompose $v$ as $v=\sum_j v_j$ with $v_j\in V_j$.

Using the symplectization defined by $\mu_\tau$ one gets easily
$$E_h(c^s_v)  =\left\{
\begin{array}{cl}
+\infty&{\rm if} \exists j\ {\rm such\ that}\ s_j>0 \ {\rm and}\ v_j\ne
0\\
\langle \tau,s\rangle -\langle \mu_\tau(v), -is\rangle&{\rm otherwise}.
\end{array}
\right.
$$

Therefore, if $E_h(c^s_v)<\infty$, one must have 
$$v_j\ne 0\Rightarrow s_j\leq 0\ .$$
But this clearly implies  that  $c^s_v(t)=\sum\limits_j e^{ts_j} v_j$ tends
to
$\sum\limits_{j,s_j=0} v_j$ as $t\ra\infty$.
\qed

\paragraph{The map $\lambda$ on the cone $H(G)$.} Fix a representative
$(K,g,\mu)\in\sigma$. We recall (see section 1.2)  that to every pair
$(s,x)\in i\kg\times X$  we associated the map
$$\lambda_x^s: \R_{\geq 0}\ra\R\ ,\ \lambda_x^s(t):=\mu^{-is}(c^s_x(t))\ ,
$$
where   $\mu^\xi:=\langle
\mu,\xi\rangle:X\ra\R$ and $c^s_x(t):=\exp(ts) x$.  If we choose an equivalent
triple
\begin{equation}\label{equivariance}
(K',g',\mu')=(\Ad_\gamma(K), (\gamma^{-1})^*
g,
\ad_{\gamma^{-1}}^t\circ \mu\circ
\gamma^{-1})
\end{equation}
 the corresponding maps $\lambda'$ are given by
\begin{equation}\label{lambda'}
(\lambda')^s_x(t)=\lambda^{\ad_{\gamma^{-1}}(s)}_{\gamma^{-1} x}(t)\ . 
\end{equation}
Moreover, using the equivariance property of the moment map with respect to
the $K$-action, one gets

\begin{equation}\label{K-eq}
\lambda^{\ad_{k }(s)}_{k x}(t)=\lambda_x^s(t) \ \ \forall k\in K\ \forall  t\in
\R 
\ .
\end{equation}

We recall that   $\lambda^s_x$ is increasing and that we   put
\begin{equation}\label{DefLambda} 
\lambda^s(x):=\lim_{t\ra\infty}\lambda^s_x(t)\in\R\cup\{\infty\}\ .
\end{equation}
so that $\lambda^s(x)=\lambda^s_x(0)+E_g(c^s_x)$ (see section 1. 2).

\begin{lm}\label{WellDef}  Suppose that $\sigma$ is energy complete. The map
$$x\mapsto \lambda^s_x:X\ra\R$$
  does not depend on the  
representative $(K,g,\mu)\in\sigma$ with $s\in i\kg$.
\end{lm}
\pf  Let $(K,g,\mu)\in\sigma$ such that $s\in i\kg$ and let $\gamma\in G$
such that $\ad_{\gamma{-1}}(s)\in i\kg$. We consider the representative
$(K',g',\mu')=(\Ad_\gamma(K), (\gamma^{-1})^* g,
\ad_{\gamma^{-1}}^t\circ \mu\circ
\gamma^{-1})$ of $\sigma$. taking into account the equivariance formula 
(\ref{equivariance}) we have to show that
\begin{equation}\label{ToShow}
\lim_{t\ra\infty}\lambda^s_x(t)=
\lim_{t\ra\infty}\lambda^{\ad_{\gamma^{-1}}(s)}_{ \gamma^{-1} x }(t)  
\end{equation}
Decompose $\gamma^{-1}$ as
$\gamma^{-1}=kh$, where $k\in K$,  $h\in\exp(i\kg)$. By   Lemma \ref{Dec},
$\ad_h(s)=s$, hence, by (\ref{K-eq}), we get
$$\lambda^{\ad_{\gamma^{-1}}(s)}_{\gamma^{-1}
x}(t)=\lambda^{\ad_{k}(s)}_{kh x}(t)=\lambda ^s_{hx}(t)\ .
$$
Therefore,   our claim (\ref{ToShow}) reduces to the  formula 
$$\lim_{t\ra\infty}\lambda^s_x(t)=
\lim_{t\ra\infty}\lambda_{hx}^s(t)\ ,
$$
for   $h\in\exp(i\kg)$ with $\ad_h(s)=s$.

Suppose first that $E_g(c^s_x)<\infty$.  Since $\alpha$ is energy complete,
  the limit
$$l=\lim_{t\rightarrow\infty} e^{ts}x\in X
$$
exists. Choose $\chi\in i\kg$ such that $h=e^\chi$. We obtain
$$\lambda^s_x(t)-\lambda_{hx}^s(t)=\mu^{-is}(e^{ts}x)-\mu^{-is}(e^{ts}hx)=
\mu^{-is}(e^{ts}x)-\mu^{-is}(he^{ts}x)=$$
$$=
-\int_{\left.c_{e^{ts}x}^\chi\right|_{[0,1]} } d(\mu^{-is})\ ,
$$
because $\left.c^\chi_{e^{ts}x}\right|_{[0,1]}$ is a curve  joining
$e^{ts} x$ to $h e^{ts} x$.

But $d\mu^{-is}(\cdot)=\omega_g((-is)^\#,\cdot)=g(s^\#,\cdot)$, hence
$$\lim_{t\rightarrow\infty}(\lambda^s_x(t)-\lambda_{hx}^s(t))=-\int
_{\left.c_{l}^\chi\right|_{[0,1]}} g(s^\#,d c_{l}^\chi)\ .
$$
Since $[\chi,s]=0$, one has
$c^\chi_l(\tau)=e^{\tau\chi}l=\lim_{t\ra\infty} e^{\tau\chi}
e^{ts}x=\lim_{t\ra\infty} e^{ts}(e^{\tau\chi}x)$, so $c^\chi_l(\tau)$ is
a fixed point of the local flow  associated with the vector field  $s^\#$.
Therefore 
$s^\#$ vanishes identically along the curve $c_{l}^\chi$, and this completes the
proof in the case $E_g(c^s_x)<\infty$.  We argue similarly when
$E_g(c^s_{hx})<\infty$; we just replace $x$ by $hx$ and $h$ by $h^{-1}$.

If finally $E_g(c^s_x)=E_g(c^s_{hx})=\infty$, we get by (\ref{energy}) that
$\lambda^s(x)=\lambda^s(hx)=\infty$.
\qed

The previous lemma plays a crucial role in developing our stability
concept: it allows us to associate to any energy complete
symplectization
$\sigma$ a well defined map
$$\lambda: H(G)\times X\ra\R\cup\{\infty\}\ ,\ (s,x)\mapsto\lambda^s(x)
$$
which is  {\it intrinsically} associated with $\sigma$ and has the following
important properties.
\begin{pr}\label{properties} Suppose that $\sigma$ is energy complete. The map
$\lambda$ introduced above satisfies the following properties:
\begin{enumerate}
\item homogeneity: $\lambda^{ts}(x)=t\lambda^s(x)$ for any $t\in \R_{\geq
0}$.
\item  equivariance:
$\lambda^s(x)=\lambda^{\ad_{\gamma }(s)}(\gamma x)$ for all $s\in
H(G)$, $\gamma\in G$.
\item  parabolic  invariance: $\lambda^s(x)=\lambda^s(h x)$ for every $h\in
G(s)$.
\item   $\sim$ invariance: $\lambda^s(x)=\lambda^\sigma(x)$ if $s\sim\sigma$.
\item  semicontinuity: if $(x_n,s_n)_n\ra (x,s)$ then $\lambda^s(x)\leq
\liminf\limits_{n\ra\infty}\lambda^{s_n}(x_n)$.
\end{enumerate}
\end{pr}
\pf  The first statement follows directly from  the
definition.   The second property follows from Lemma \ref{WellDef} and
(\ref{lambda'}). To prove the third, use (\ref{SemProd}) and put
$h=h_0 h'$, where $h_0$, $h'\in G$ with
$\ad_{h_0} s=s$ and
$\lim_{t\ra\infty} e^{ts} h' e^{-ts}=e$. We get by 1.
$$\lambda^s(hx)=\lambda^s(h_0h'x)=\lambda^{\ad_{h_0^{-1}}(s)}(h'x)=
\lambda^s(h'x)\ ,
$$
so it  remains to show that $\lambda^s(h'x)=\lambda^s(x)$. 

Suppose that $\lambda^s(x)<\infty$. By energy completeness, this implies that
the limit $l=\lim_{t\ra\infty} e^{ts}x$ exists in $X$.  Therefore
$$\lim_{t\ra\infty} e^{ts} h'x=\lim_{t\ra\infty} e^{ts} h' e^{-ts} e^{ts}x =l\ ,
$$
because $e^{ts} h' e^{-ts}\ra e$. Therefore
$$\lambda^s(h'x)=\lim_{t\ra \infty} \mu^{-is}(e^{ts}
h'x)=\mu^{-is}(l)=\lambda^s(x)\ .
$$
In the case $\lambda^s(h'x)<\infty$ we argue similarly. This completes the proof
of the third property.

The fourth follows immediately from 2. and 3. 

For the fifth, note first that we may suppose that there exists a maximal
compact subgroup $K$ of $G$  such $s_n\in i\kg$. If our conclusion was false,
there would exist $\varepsilon>0$ and a subsequence
$(x_{n_m},s_{n_m})_m$ of $(x_n,s_n)_n$ such that the limit $\lim_{m\ra\infty}
\lambda^{s_{n_m}}(x_{n_m})$ exists, is finite and
$$\lambda^s(x)\geq
\lim_{m\ra\infty}
\lambda^{s_{n_m}}(x_{n_m})+\varepsilon\ .$$
   Fix a representative
$\rho=(K,g,\mu)\in\sigma$ and choose
$t$ sufficiently large such that, with respect to $\rho$, one has
$\lambda^s_x(t)\geq
\lim_{m\ra\infty}\lambda^{s_{n_m}}(x_{n_m})+\frac{1}{2}\varepsilon$. But,
since
$$\lambda^{s_{n_m}}(x_{n_m})\geq  \lambda^{s_{n_m}}_{x_{n_m}}(t)$$
and $(s,x)\mapsto \lambda^s_x(t)$ is continuous on $i\kg\times X$, one gets
$$\lim_{m\ra\infty}\lambda^{s_{n_m}}(x_{n_m})\geq
\lim_{m\ra\infty}\lambda^{s_{n_m}}_{x_{n_m}}(t)=\lambda^s_{x}(t) \geq
\lim_{m\ra\infty}\lambda^{s_{n_m}}(x_{n_m})+\frac{1}{2}\varepsilon\ .
$$
which is a contradiction.
\qed
\subsection{Analytic stability, polystability, semistability}

Let again $\alpha$ be an action of a reductive group $G$ on a
complex manifold $X$. Fix a
symplectization $\sigma$ of $\alpha$ and let $\lambda:H(G)\times X\ra
\R\cup\{\infty\}$ be the associated map.

Let $G$ be a reductive group. A subalgebra $\g'\subset\g$ will be called
a  {\it reductive subalgebra} if it has the form $\g'={\kg'}^\C$,
where $\kg'$
is the Lie algebra of a compact subgroup of $G$.  Equivalently, one can
require instead that $\g'=(\kg\cap\g')^\C$ for a maximal compact subgroup $K$
of $G$. Note that, if
$\kg'$ is the Lie algebra of a compact Lie subgroup of $G$, then the minimal
complex subspace
$\tilde \kg'$ of $\g$ which contains $\kg'$ can be identified with the
complexification ${\kg'}^\C$ via the canonical map ${\kg'}^\C\ra \tilde \kg'$.

\begin{dt}\label{AnStab} A
point $x\in X$ will be called
\begin{enumerate}
\item analytically $\sigma$-\ub{semistable}  if
$\lambda^s(x)\geq 0
$
for all $s \in H(G)$.
\item analytically $\sigma$-\ub{stable} if it is semistable
   and
$\lambda^s(f)>0$ for  $s\in H(G)\setminus\{0\}$.
\item analytically $\sigma$-\ub{polystable} if it is semistable,  $\g_x$ is
a reductive subalgebra\footnote{We will see that analytic polystability
implies symplectic polystability, hence the stabilizer $G_x$ of a
polystable point $x$ is reductive by Proposition \ref{RedStab}. We preferred to
require only the reductivity of   $\g_x$ in our definition in order to have a
purely infinitesimal condition.},
  and $\lambda^s(x)> 0$ if  $s$  is not
equivalent to an element of $\g_x$.
\end{enumerate}
\end{dt}
\begin{re} Let $x\in X$. If  $x$ is analytically
$\sigma$-semistable, then
$s\mapsto\lambda^s(x)$ must vanish   on $H(G)\cap \g_x$.
\end{re}
\pf  For every $(K,g,\mu)\in\sigma$ and $s\in i\kg\cap\g_x$ one
has $\lambda^s(x)=\mu^{-is}(x)$. This shows that 
$\lambda(x)^{ts}=t\lambda(x)^{s}$ for every $s\in\g_x\cap H(G)$ and
$t\in\R$. The semistability condition implies that $\lambda^s(x)\geq 0$  and
$-\lambda^s(x)=\lambda^{-s}(x)\geq 0$ for every $s\in  H(G)\cap \g_x$.
\qed

The following proposition shows that it suffices to check the (semi, poly-)
stability conditions   for vectors $s\in i\kg$, where $\kg$
is the Lie algebra of a fixed maximal compact subgroup. However, one should be
very careful in the polystable case.

\begin{pr}\label{K-criterion} Suppose that $\sigma$ is energy complete.
\begin{enumerate}
\item  The analytical $\sigma$-semistability (stability,
polystability) condition for
$x\in X$ depends only on the complex orbit $Gx$ of $x$.
\item Choose any maximal compact subgroup $K\subset G$. A point $x\in X$ is

\begin{enumerate}
\item analytically $\sigma$-semistable, if and only if
$\lambda^s(x)\geq0$ for every $s\in i\kg$.
\item analytically $\sigma$-stable,
  if and only if $\lambda^s(x)\geq 0$ for any $s\in i\kg$ and
$\lambda^s(x)>0$ when $s\in i\kg\setminus\{0\}$.
\item analytically $\sigma$-polystable,   if and only if $\lambda^s(x)\geq 0$ for
any $s\in i\kg$, $\g_x$ is a reductive Lie algebra, and
$\lambda^s(x)>0$ for any
$s\in i\kg$ which is not equivalent to an element of $\g_x$.
\end{enumerate}
\item Suppose that   $\g_x$   is a reductive Lie algebra  and let  $\kg$ be
the Lie algebra of a maximal compact subgroup of $G$ such that
$\g_x=(\kg\cap\g_x)^\C$. 
   Then
$x$ is analytically $\sigma$-polystable if and only if
$\lambda^s(x)\geq 0$ for any
$s\in i\kg$ and $\lambda^s(x)>0$ when $s\in i\kg\setminus i\kg_x$.
\end{enumerate}
\end{pr}
\pf The first statement follows from the equivariance property 2. in
Proposition  \ref{properties}.  The second follows   from the fact that $i\kg$
is a complete system of representatives for the relation $\sim$ (Corollary
\ref{Ret}) and the same equivariance property.

For the third statement, note that,     when
$\g_x=(\kg\cap\g_x)^\C$,   any element of
$i\kg$ which is equivalent in
$H(G)$ to  an element of $H(G)\cap \g_x$ must belong to $i\kg\cap \g_x$. This
follows from Corollary \ref{embedding}, applied to the inclusion of reductive
groups
$(K\cap G_x)^\C\subset G$.
\qed

\section{Comparison Theorems. Hilbert criterion in   K\"ahlerian geometry}

Here is our first comparison result.
\begin{thry}\label{FirstComp} Suppose that $\sigma$ is energy complete. A point
$x\in X$ is symplectically  $\sigma$-stable (polystable)  if and only if it
is analytically
$\sigma$-stable (polystable).
\end{thry}
\pf  The stable case follows immediately from Proposition \ref{K-criterion} and
Proposition \ref{NumCrit}.

{\it "symplectically  polystable" $\Rightarrow$ "analytically 
polystable"}:  

Fix a
representative $(K,g,\mu)\in\sigma$. 
 Let $x \in X$ be a symplectically polystable point.   Let
$x_0\in Gx\cap \mu^{-1}(0)$. By Proposition
\ref{K-criterion}, it suffices to prove that $x_0$ is analytically  polystable. But,
from the proof of Proposition \ref{NumCrit}, we know that
$\g_{x_0}=\kg_{x_0}^\C$  hence, by Proposition \ref{K-criterion}, it suffices to
prove that
$\lambda^s(x_0)\geq 0$ for all $s\in i\kg$ and $\lambda^s(x_0)> 0$ for
all $s\in\kg\setminus\g_{x_0}$. But $\lambda^s(x_0)\geq
\lambda^s_{x_0}(0)=0$, because
$\mu(x_0)=0$. Moreover, for   $s\not\in \g_{x_0}$ one has
$s^\#(x_0)\ne 0$, which
implies that
$$\lambda^s(x_0)>\lambda^s_{x_0}(0)=0\ .$$

"analytically  polystable" $\Rightarrow$ "symplectically  polystable":

Since $x$ is analytically polystable, the  subalgebra $\g_x$ is reductive.
Choose a maximal compact subgroup $K$ of $G$ such that
$\g_x=(\g_x\cap\kg)^\C$.   By Proposition
\ref{K-criterion}, 3. we know that   
$\lambda^s(x)\geq 0$ for all $s\in i\kg$ and  $\lambda^s(x)> 0$ for
$s\in i\kg\setminus i\kg_x$.  It suffices to apply Proposition \ref{NumCrit}
using 
  a representative $\rho$ of the symplectization
$\sigma$ with first component  $K$ and to   remember that symplectic
polystability   depends only on the fixed symplectization. 
\qed
\vspace{3mm}\\
{\bf Remark:} The assumption that $\sigma$ is energy 
 complete plays implicitly an important role in this proof. It gives us the
flexibility to choose the maximal compact subgroup of $G$ in a convenient way.
\begin{lm}\label{orbit}  Suppose that  $\sigma$ is energy complete. Let
$x_0\in X$ be an analytically $\sigma$-semistable point, and let $x\in X$ such
that $x_0\in\overline{Gx}$. Then $x$ is analytically $\sigma$-semistable.
\end{lm}
\pf
Suppose that $x$ was not semistable, hence there exists $s\in H(G)$ such that
$\lambda^s(x)<0$. Let $(g_n)_n$ be a sequence in $G$ such that $(g_n x)_n$
converges to $x_0$. Using the equivariance property of the function
$\lambda$ (Proposition \ref{properties}), one gets
$$\lambda^{\ad_{g_n}(s)}(g_n x)=\lambda^s(x)<0\ .
$$
Let $K$  be a maximal compact subgroup of $G$ such that $s\in i\kg$ and let
$s_n\in i\kg$ such that $s_n\sim \ad_{g_n}(s)$ (see Proposition \ref{Ret}).
By Proposition \ref{CompactOrbit}, one can find a subsequence $(s_{n_m})_m$ of
$(s_n)_n$ which converges to a vector, say $s_0$, in the orbit  $\ad_{K}(s)$. By
the semicontinuity property and the $\sim$ invariance property in Proposition
\ref{properties}, we obtain
$$\lambda^{s_0}(x_0)\leq \liminf \lambda^{s_{n_m}}(g_{n_m}x)=\liminf
\lambda^{\ad_{g_{n_m}}(s)}(g_{n_m} x)=\lambda^s(x)<0\ ,
$$
which contradicts the analytic semistability of $x_0$.
\qed

 Our next comparison theorem is more delicate:

\begin{thry}\label{semistability}
Let $(X,g)$ be a complex manifold, $\alpha : G \times X \to X$ a complex
reductive Lie group action and let $\sigma=[K,g,\mu]$ be an energy complete
symplectization of this action.   Then,
for any point
$x\in X$ the following properties are equivalent:
\begin{enumerate}
\item   $x$ is symplectically $\sigma $-semistable;
\item  $x$ is analytically  $\sigma $-semistable.
\item $\lambda^s(x)\geq 0$ for all $s\in iz_{\kg}(\kg_x)$.
\item There exist  $s_m,\
s_0\in iz_{\kg}(\kg_x)$ such that 
\begin{enumerate}
\item $\lambda^{s_m}(x)=0$,
\item $[s_m,s_0]=0$, 
\item  the limit $y=\lim_{t\ra\infty}(\exp(its_m)x)$
exists in $X$, 
\item  $\mu(\exp(s_0)y)=0$.
\end{enumerate}
\end{enumerate}
\end{thry}
\pf  \\
1.$\Rightarrow$ 2.

We fix a representative
$(K,g,\mu)\in\sigma$. Since $x$ is symplectically $\sigma$-semistable,  there
exists a point $x_0\in \overline{Gx}\cap \mu^{-1}(x_0)$.  By Theorem
\ref{FirstComp}, $x_0$ is analytically $\sigma$-polystable, in particular it is
analytically $\sigma$-semistable. It follows that $x$ is also analytically
$\sigma$-semistable, by Lemma \ref{orbit}.
\\ \\
2.$\Rightarrow$ 3. Obvious
\\ \\
3. $\Rightarrow$ 4.

 Fix a representative $(K,g,\mu)\in\sigma$.

Since 3. holds, one of the following two
possibilities must occur:
\begin{enumerate} 
\item  [A.] For all $s\in iz_{\kg}(\kg_x)\setminus iz(\kg_x)$ it holds
$\lambda^s(x)>0$.
\item  [B.] There exists $s\in iz_{\kg}(\kg_x)\setminus iz(\kg_x)$ such that
$\lambda^s(x)=0$.
\end{enumerate}

In the first case,   Lemma \ref{lemma} gives an element
$s_0$ in the orthogonal complement $i[z(\kg_x)^{\bot_{z_{\kg}(\kg_x)}}]$  of
$iz(\kg_x)$ in $iz_{\kg}(\kg_x)$   such that $\mu(\exp(s_0)x)=0$.
Therefore, taking
$s_m=0$, the claim 4. holds   it this case.

 In the second case we know that there exist vectors $s\in
iz_{\kg}(\kg_x)\setminus iz(\kg_x)$  such that $\lambda^s(x)=0$. For any
such $s$, the limit $x_s=\lim_{t\ra\infty} (\exp(ts)x)$ exists, because
$\sigma$ is energy complete. We choose $s_m \in iz_{\kg}(\kg_x)\setminus
iz(\kg_x)$ with
$\lambda^{s_m}(x)=0$   which maximizes
$\rk(K_{x_s})$, i. e. such that
\begin{equation}\label{maximize}
\rk(K_{x_{s_m}})\geq \rk(K_{x_s})\ \forall  s  \in
iz_{\kg}(\kg_x)\setminus iz(\kg_x) \hbox{ with }\lambda^{s}(x)=0\ .
\end{equation}
 We remind that the rank of a compact Lie
group is the dimension of a maximal torus of it.  Set
$y:=x_{s_m}$.

Note that $s_m\in\g_y$, so $is_m\in \kg_y$.   Moreover,   since $\exp(ts)\in
Z_K(\kg_x)$, it follows that  $\kg_x\subset \kg_{\exp(ts)x}$ for all
$t\in\R$, hence, making $t\ra\infty$, one gets $\kg_x\subset\kg_y$. 
Therefore,
\begin{equation}\label{k_y}
\{is_m\}\cup k_x\subset \kg_y\ .
\end{equation}

{\bf Claim:} $\lambda^s(y)\geq 0$, for all $s\in iz_{\kg}(\kg_y)$.

Indeed, if there   existed $\tau\in iz_{\kg}(\kg_y)$ with
$\lambda^\tau(y)<0$, one would get  immediately a contradiction in the
following way:  

By energy completeness, the limit $z=\lim_{t\ra\infty} (\exp(t\tau)y)$
exists. Consider any compact torus
$T$ of  
$K_y$ which contains the real one parameter subgroups generated by $is_m$
and
$i\tau$. Such a torus exists, because  $\tau\in iz_{\kg}(\kg_y)$ and $s_m\in
i\kg_y$ by (\ref{k_y}) so
$[\tau,s_m]=0$.\footnote{Note that, in general, in  a connected
compact group $L$ the stabilizer $Z_L(u)$ of any element $u\in\lg$ is also
connected. Moreover, $Z_L(u)$ is just the union of all maximal tori of $L$
whose Lie algebra contain the vector $u$.} so
$[\tau,s_m]=0$. As in the proof of Proposition  \ref{EnCompl}, one  can linearize
the induced
$T$-action on an open ${T}^\C$-invariant Stein neighborhood of $z$ (which
will necessarily contain $y$ and $x$).  Using such a linearization and the
standard eigenspace decomposition associated with a torus action on a vector
space, one can  see easily that
$$\lim_{t\ra\infty}\exp (t(s_m+  \varepsilon \tau)x)=z$$
 for all sufficiently
small $\varepsilon>0$.  This implies that
\begin{equation}\label{sum}
\lambda^{s_m+\varepsilon\tau}(x)=\mu^{-i(s_m+\varepsilon\tau)}(z)=
\mu^{-i s_m}(z) +\varepsilon \mu^{-i \tau}(z)=\mu^{-i s_m}(z) +\varepsilon
\lambda^\tau(y)\ .
\end{equation}

The point $y=x_{s_m}$ is the limit for $t\ra\infty$ of an integral curve  of
the vector field $s^\#_m$, hence   $y$ is  invariant under the
local one-parameter transformation group generated by  
$s^\#_m$. Since $s_m$ and
$\tau$ commute,  all the points of the curve $c^\tau_{y}$ will be
fixed under this local one-parameter transformation group. Therefore 
\begin{equation}\label{s}
s^\#_m|_{\im (c^\tau_{y})}\equiv 0\ .
\end{equation}

On the other hand, by the properties of a moment map, one has
$$d\mu^{-is_m}=\iota_{-is^\#_m}\ \omega_g\ ,
$$
which, together with (\ref{s}), shows that the real function $\mu^{-is_m}$ is
constant on the curve $c^\tau_{y}$ which joins $y$ to $z$. On the other
hand,  taking into account the way in which $s_m$ was chosen,   one has 
$$\lambda^{s_m}(x)=
\mu^{-i s_m}(y)=0\ .$$
 Therefore  $\mu^{-is_m}(z)=0$  hence, by (\ref{sum}),
$\lambda^{s_m+\varepsilon\tau}(x)<0$, which contradicts the analytic 
semistability of $x$. This proves the claim. \\ \\

We can distinguish now again the following two cases:
\begin{enumerate} 
 \item [a.]  For all $\tau\in i z_{\kg}(\kg_y)\setminus i z(\kg_y)$
it holds $\lambda^\tau(y)>0$.
\item [b.] There exists $\tau\in i z_{\kg}(\kg_y)\setminus i z(\kg_y)$
such that $\lambda^\tau(y)=0$.
\end{enumerate}

In the case a., using again Lemma
\ref{lemma} we get a vector $s_0\in i z(\kg_y)^{\bot_{z_\kg(\kg_y)}}$ such
that
\begin{equation}\label{weak} \mu(\exp(s_0)y)=0\ .
\end{equation}  

Since $s_m\in i\kg_x\subset i\kg_y$ and $s_0\in iz_\kg(\kg_y)$, we have
$[s_m,s_0]=0$, hence  4. is proved in  case a.  

We will
show now that in, fact  case b. cannot occur, because it would contradict the
maximizing property (\ref{maximize}). Indeed, if b. held, let  
$\tau\in i z_{\kg}(\kg_y)\setminus i z({\kg_y})$ such that
$\lambda^\tau(y)=0$.  By energy completeness, the curve $c^\tau_y$ has a limit
$z$ for $t\ra \infty$  and, similarly to (\ref{k_y}), we will have
\begin{equation}\label{k_z}
\{i\tau\}\cup\kg_y\subset\kg_z.
\end{equation}

 The same method as in the proof of
the claim above, shows that
$$\lim(\exp(t(s_m+\varepsilon\tau)x)=z$$
 for all sufficiently small
$\varepsilon>0$ and that (\ref{sum}) holds. This time we obtain
$\lambda^{s_m+\varepsilon\tau}(x)=0$ because $\lambda^\tau(y)=0$. Let $T_m$
be a maximal compact torus of $K_y$. The stabilizer
$K_z$ of $z$ contains both $T_m$ (because   $K_y^e\subset K_z$ by
(\ref{k_z})) and the real one parameter subgroup generated by
$i\tau$. The elements of this one parameter subgroup commute
with the elements of $T_m$ because $\tau\in z_{\kg}(\kg_y)$.

We know that $\tau \in i z_{\kg}(\kg_y)\setminus i z({\kg_y})=i
z_{\kg}(\kg_y)\setminus \kg_y$, hence $\tau\not\in\kg_y$, hence
$i\tau\not\in\tg_m$.  it follows that $T_m$ and this real
one parameter subgroup  generate a torus in $K_z$ which has a larger
dimension that $\dim(T_m)$. This contradicts the maximizing property
(\ref{maximize}).

4.$\Rightarrow$ 1. Using  4. we get
$$\exp(s_0)\lim_{t\ra\infty}(\exp(ts_m)x)\in\overline{Gx}\cap\mu^{-1}(0)\ ,$$
hence
$x$ is symplectically $\sigma$-semistable.
\qed
\\ \\

We can prove now the following important semistability criterion.  Note that
this result is obvious in the case of a compact manifold $X$, whereas for non
energy complete symplectizations it  is in general false.
\begin{thry}\label{criterion} 
Let $(X,g)$ be a complex manifold, $\alpha : G \times X \to X$ a complex
reductive Lie group action and let $\sigma $ be an energy complete
symplectization of this action.   Then,
for any point
$x\in X$ the following properties are equivalent:
\begin{enumerate}
\item   $x$ is symplectically $\sigma $-semistable;
\item  $x$ is analytically  $\sigma $-semistable.
\item $\inf\limits_{g\in G}  \|\mu(gx)\| =0$, where the norm is computed with
respect to any
$\ad$-invariant inner product on $\kg$.

\end{enumerate}
\end{thry}
\pf We know that 1. and 2. are equivalent.  The implication 1.$\Rightarrow$ 3.
is obvious, so it suffices to show that 3.$\Rightarrow$2.

Suppose that $x$ was not analytically semistable, and fix a triple
$(K,g,\mu)\in\sigma$.  It would follow that there exists $s\in i\kg $ such that
$\lambda^s(x)<0$. We normalize $s$ such that $\| s\|=1$.  Since 3. holds, there
exists a sequence $(g_n)_n$ in $G$ such that $\|\mu(g_nx)\|\ra 0$.  By the
equivariance  property in Proposition \ref{properties}, we get 
$\lambda^{\ad_{g_n}(s)}(g_n x)=\lambda^s(x)$. Using Corollary \ref{Ret},   we can
find $s_n\in i\kg$ such that $s_n\sim  \ad_{g_n}(s)$. 

Using the $\sim$ invariance property, we get for all $n\in\N$.
\begin{equation}\label{negative}
\mu^{-is_n}(g_n
x)=\lambda_{g_n x}^{s_n}(0)\leq\lim_{t\ra\infty}\lambda_{g_n x}^{s_n}(t)=
\lambda^{s_n}({g_n x})=\lambda^s(x)<0\ .
\end{equation}
The point is now that, by
Proposition \ref{CompactOrbit}, $s_n\in \ad_{K}(s)$, hence $\|s_n\|=\|s\|=1$.
Thus,  
$$|\mu^{-is_n}(g_n x)|= |\langle \mu(g_n x),-is_n\rangle|\leq \|\mu(g_n x)\|\ra
0\ ,
$$
and this obviously contradicts (\ref{negative}).
\qed

Our last semistability criterion is the following:

\begin{thry}\label{lastcriterion} 
Let $(X,g)$ be a complex manifold, $\alpha : G \times X \to X$ a complex
reductive Lie group action and let $\sigma=[K,g,\mu]$ be an energy complete
symplectization of this action.   Then,
for any point
$x\in X$ the following properties are equivalent:
\begin{enumerate}
\item   $x$ is symplectically $\sigma $-semistable;
\item  $x$ is analytically  $\sigma $-semistable.
\item The real function $\Psi(x,\cdot)$ is bounded from below on $G$.
\end{enumerate}
\end{thry}
\pf

The implication 3.$\Rightarrow$ 2. is very easy: if there existed $s\in i\kg$ such
that $\lambda^s(x)<0$, then, taking into account that
$$\frac{d}{dt}\Psi(x,\exp(ts))=\lambda^s_x(t)\leq\lambda^s(x)\ ,
$$
we get an estimate of the form $\Psi(x,\exp(ts))\leq C_1+\lambda^s(x) t$, which shows
that $\lim_{t\ra\infty}\Psi(x,\exp(ts))=-\infty$.

We prove now 1.$\Rightarrow$ 3. 

Let $x$ be a symplectically polystable point. By the results of [HH] (see Theorem
\ref{KQuot}), the closure of the orbit
$Gx$ contains a unique polystable orbit $O$, and this orbit contains a unique
$K$-orbit $o$ on which $\mu$ vanishes.

Now we use the
existence theorem for local potentials compatible with a Hamiltonian action  
([HHL], p. 138,  [HH], p. 245). There exists a
$G$-invariant    Stein neighborhood $U$ of the $K$-orbit $o$ and  a {\it potential
for the moment map
$\mu|_U$}, i. e. a strictly plurisubharmonic
$K$-invariant function $\varphi:U\ra\R$ such that on $U$ it holds
\begin{equation}\label{varphi}
\langle\mu,a\rangle=d\varphi(J a^\#)\ ,\ \forall a\in\kg\ ,
\end{equation}
where $J\in\End(T_X)$ is the almost complex structure of the complex manifold $X$.  

We state that (\ref{varphi}) implies   the following formula for the
restriction $\Psi|_{U\times G}$ of function
$\Psi$ associated with $\mu$:
\begin{equation}\label{PsiPot}\Psi(u,g)=\varphi(gu)-\varphi(u)\ ,\ \forall u\in U,\
\forall g\in G\ .
\end{equation}
Indeed, using (\ref{varphi}), we get
$$\frac{d}{dt}(\varphi(e^{ts}u))=d\varphi(s^\#_{e^{ts}u})=d\varphi(J(-is)^\#_{e^{ts}u})=
\langle\mu(e^{ts}u)),-is\rangle=\lambda^s_{e^{ts}u}(t)\ .
$$
The other conditions in Definition \ref{Psi} are obvious.

Let $x_0\in o$. Since  $\mu(x_0)=0$,    by Remark \ref{Crit} it follows that  $e\in G$
is a critical point of the map $\Psi(x_0,\cdot)$.  Taking into account the convexity
property of the map
$\Psi$ (Remark \ref{PsiProp}) and Proposition \ref{RedStab}, it follows that  $\Psi(x_0,\cdot)$
reaches its absolute minimum in $e$. Therefore, by (\ref{PsiPot}) we get  
\begin{equation}\label{ineg}
\varphi(x_0)\leq
\varphi(y)\ ,\ \forall  y\in O\   .
\end{equation}

By Theorem \ref{semistability}, for any  $\xi\in Gx$ there exist $s_m^\xi\in i\kg$  and
$y^\xi\in O$ such that %
\begin{equation}\label{form}\lambda^{s_m^\xi}(\xi)=0\ ,\ 
\lim_{t\ra\infty} \exp(ts_m^\xi)=y^\xi\ .
\end{equation}

Since $U$ is open,  $G$-invariant  and it contains the $G$-orbit $O$, it follows that $Gx\subset
U$, so we can write
$$\Psi(\xi,\exp(ts_m^\xi))=\varphi(\exp(ts_m^\xi))-\varphi(\xi) \ra \varphi(y^\xi)-\varphi(\xi)\
{\rm as}\ t\ra\infty.
$$
The first relation in (\ref{form}) shows that  $\lambda^{s_m^\xi}_{\xi}(t)\leq 0$ forall 
$t\in\R_{\geq 0}$, hence the real function  $t\mapsto \Psi(\xi,\exp(ts_m^\xi))$ is decreasing,
so
$$\varphi(y^\xi)-\varphi(\xi)\leq \Psi(\xi,e)=0\ .
$$
Combining with (\ref{ineg}), we get
$$\varphi(\xi)\geq \varphi(y^\xi)\geq \varphi(x_0)\ ,
$$
hence $\Psi(x,g)=\varphi(gx)-\varphi(x)\geq \varphi(x_0)-\varphi(x)$.
\qed

\vspace{3mm}
{\small
Author's address: \vspace{2mm}\\
LATP, CMI,   Universit\'e de Provence,  39  Rue F. Joliot-Curie, 13453 Marseille
Cedex 13, France,  e-mail: teleman@cmi.univ-mrs.fr 
    ,  and\\ Faculty of Mathematics, University of Bucharest, Bucharest,
Romania }
 
\end{document}